\documentclass[preprint,11pt]{elsarticle}

\usepackage[latin1]{inputenc}
\usepackage{amssymb}
\usepackage{amsmath} 
\usepackage{amsfonts}
\usepackage{rawfonts}
\usepackage{oldlfont}
\usepackage{amsthm}
\usepackage{enumerate}
\usepackage{graphicx}
\usepackage[dvips]{epsfig}
\usepackage{epsfig}
\newtheorem{lemma}{Lemma}
\newtheorem{lemma-a}{Lemma}[subsection]
\newtheorem{th-a}{Theorem}[subsection]

\newtheorem{definition}{Definition}

\newtheorem{cor-a}{Corollary}[subsection]

\def\R{{\mathbb{R}}}
\def\Ll{{\mathcal{L}}}
\def\P{{{\mathrm{P}}}}
\def\E{{{\mathbb{E}}}}
\def\X{{\mathcal{X}}}
\def\ni{\noindent}
\def\argmin{{\mathop{{\mathrm{argmin}}}}}
\def\argmax{{\mathop{{\mathrm{argmax}}}}}

\journal{Statistics and Probability Letters}

\begin{document}

\begin{frontmatter}

\title{Vapnik-Chervonenkis Dimension of Axis-Parallel Cuts}
\author[rvt]{Servane Gey}
\ead{Servane.Gey@parisdescartes.fr}
\address[rvt]{Laboratoire MAP5 - UMR 8145, Universit\'e Paris Descartes, 75270 Paris Cedex 06, France}

\begin{abstract}
Algorithms in high dimension uses axis-parallel cuts to partition $\R^d$ in order to reduce the computational time of classifiers or regressors. Evaluating the complexity of such partitions is then crucial to evaluate estimation performance.\\ In this framework, we show that the Vapnik-Chervonenkis dimension (VC dimension) of the set of half-spaces of
$\R^d$ with frontiers parallel to the axes is of the order of $\log_2{d}$.
\end{abstract}


\end{frontmatter}

\section{Introduction}

The VC dimension of a set of subsets has been introduced by Vapnik and Chervonenkis \cite{VapCher71, VapCher74} to measure its complexity. The VC dimension of a real-valued function space $\mathcal{F}$ is then the VC dimension of $\left\{\{x ; \ f(x)\geqslant 0\} ; \ f\in \mathcal{F}\right\}$. In particular, the VC dimension of sets of classifiers or regressors appears commonly in the statistical learning area when evaluating their performance.\\ 
For example, Vapnik's theory in the classification framework is now widely known (see \cite{DevGyoLug96} for instance): let $(X,Y)$ be a couple of variables taking values in $\R^d\times \{0;1\}$, and let $\Ll$ be a sample of $n$ independent replications of $(X,Y)$. If $\hat{f}$ is a classifier minimizing the average misclassification rate of $\Ll$ on a set of classifiers having finite VC dimension $V$, then, without further assumption on the distribution $\P$ of $(X,Y)$, the performance of $\hat{f}$ is evaluated as follows: 
\begin{eqnarray} \label{riskbound}
\E_{\Ll}\left[\P\left(\hat{f}(X)\neq Y\right)\right] \leqslant C_1bias^2(\hat{f})+C_2\sqrt{\frac{V}{n}},
\end{eqnarray}
where $\E_{\Ll}$ denotes the expectation with respect to the sample distribution, $bias(\hat{f})$ denotes the bias of the classifier $\hat{f}$, and $C_1$ and $C_2$ are absolute constants.\\

Functional estimates defined on partitions of $\R^d$ are often used to estimate relationships between two variables $X\in \R^d$ and $Y\in \{0;1\}$ or $Y\in \R$ (such as histograms, piecewise polynomials, or splines for example). In many cases, the VC dimension of the set of subsets used to construct the partition appears inside risk bounds when evaluating the performance of such estimators. For example, if the set used is the set of all half-spaces of $\R^d$, often its VC dimension $d+1$ has to be taken into account.\\ When $d$ is large, it is often computationally easier to construct partitions using axis-parallel cuts. For example, some theoretical developments on dyadic partitions of $\R^2$ are given in \cite{Don97, Aka11}, and the VC dimension of axis-parallel cuts appears more particularly in the results obtained on the performance  of classification and regression binary decision trees (CART) introduced by Breiman {\it et. al} \cite{Brei84} in 1984, and theoretically studied in \cite{Nob02, GeyNed05, Gey12, GeyMar12}. In particular, it is to be found in the results of \cite{GeyMar12} that the VC dimension of axis-parallel cuts is of order $\log_2 d$. 

\section{Reminder about VC Dimension}

The VC dimension of a set ${\mathcal{A}}$ of subsets of some measurable space
$\X$ is based on counting the number of intersects of ${\mathcal{A}}$ with a
finite set of fixed points in $\X$. 
%

\begin{definition}[{\bf Vapnik-Chervonenkis Dimension}] \label{VCDim}
Let ${\mathcal{A}}$ be a set of subsets of some measurable space $\X$. Then
$(x_1,\ldots,x_n)\in \X^n$ will be said to be {\it shattered} by
${\mathcal{A}}$ if all subsets of $\{x_1;\dots;x_n\}$ are covered by ${\mathcal{A}}$, that is if $\left|\left\{\{x_1,\ldots,x_n\}\cap A \ ; \ A\in {\mathcal{A}}\right\}\right|=2^n$.\\
The {\it Vapnik-Chervonenkis dimension} $VC({\mathcal{A}})$ of ${\mathcal{A}}$ is then defined as
the maximal integer $n$ such that there exists $n$ points in $\X$ shattered by
${\mathcal{A}}$, {\it i.e.} $$VC({\mathcal{A}})=\max\left\{n \ ; \
 \max_{(x_1,\ldots,x_n)\in \X^n}\left|\left\{\{x_1,\ldots,x_n\}\cap A \ ; \
     A\in {\mathcal{A}}\right\}\right| = 2^n\right\}.$$
If no such $n$ exists, then $VC({\mathcal{A}})=+\infty$.
\end{definition}
\ni Thus, it is easily seen that the larger $VC({\mathcal{A}})$, the more
complex ${\mathcal{A}}$.\\

\ni For example, if ${\mathcal{A}}=\left\{]-\infty;x] \ ; \ x\in \R\right\}$,
$VC({\mathcal{A}})=1$; or if ${\mathcal{A}}$ is the set of all half-spaces in
$\R^d$, then $VC({\mathcal{A}})=d+1$.\\ Since axis-parallel cuts is a subset
of the set of all half-spaces in $\R^d$, it could be natural to think that its
VC dimension is of order $d$. Actually, it is shown in what follows that it is
of order $\log_2{d}$

%

\section{VC Dimension of axis-parallel cuts}

We give a formula to compute the VC dimension of axis-parallel cuts in $\R^d$. Since the obtained formula is not always easy to handle, an approximation is also given.
  
\begin{lemma} \label{lemma:vcdim}
Let $${\mathcal{A}}_d=\left\{\{x\in {\mathbb{R}}^d \ ; \ x^i\leq a\} ; \ i=1,\ldots,d \ , \ a\in
  {\mathbb{R}}\right\}.$$ Then $$VC({\mathcal{A}}_d)=\max\left\{n \ ; \
\binom{n}{\lfloor n/2\rfloor}\leq d\right\},$$ where $\lfloor n/2\rfloor$ denotes the integer part of $n/2$.\\
Furthermore, the following approximation of $VC({\mathcal{A}}_d)$ is available
for all $d\geqslant 3$:
$$ \log_2{d}+\frac{\log_2{\pi}-1}{2}\leqslant VC({\mathcal{A}}_d) \leqslant \frac{3}{2}\log_2{d}+0.63.$$
\end{lemma}
\ni {\bf Remark}: A simple calculation gives $VC({\mathcal{A}}_d)=d$ for $d
\leqslant 3$.\\

\begin{figure}[ht]
\begin{center}
\includegraphics[scale=.5]{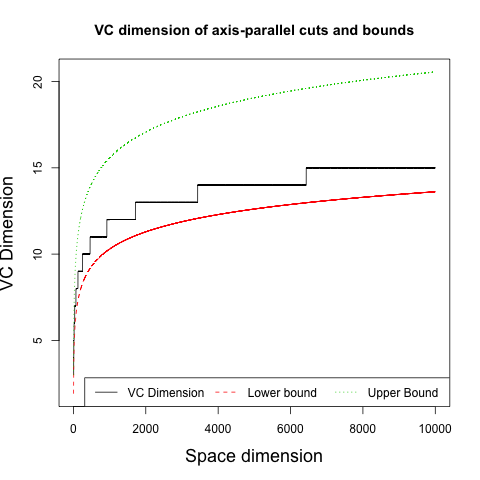}
\caption{\footnotesize{$VC({\mathcal{A}}_d)$ and Stirling's bounds with respect to the space dimension $d$.}} \label{fig:bounds}
\end{center}
\end{figure}

\ni Figure \ref{fig:bounds} shows that $VC({\mathcal{A}}_d)$ is a piecewise constant function of the space dimension $d$, which increases at a rate of order $\log_2{d}$. It also shows that the lower bound of Lemma \ref{lemma:vcdim} is conveniently sharp; the upper bound is sharp for $d$ small, and then grows farther apart from $VC({\mathcal{A}}_d)$. The bounds are obtained thanks to the Stirling's formula, which is really sharp. Actually, an approximation factor depending on $d$ has to be calibrated, leading to the observed behavior when $d$ grows.

\begin{proof}
Let $n\geqslant 1$ and $(x_1,\ldots,x_n)$ be $n$ points in ${\mathbb{R}}^d$. The idea is that, if there exists $p\leq n$ such that there is more
than $d+1$ subsets of $\{x_1,\ldots,x_n\}$ having $p$ elements, then ${\mathcal{A}}_d$ will miss at least $\binom{n}{p}-d$ subsets: suppose that $n$ is such that $\displaystyle{\binom{n}{\lfloor n/2\rfloor}> d}$. This
means that there are at least $d+1$ subsets of $\{x_1,\ldots,x_n\}$ of size
$\lfloor n/2\rfloor$. For each coordinate $i=1,\ldots,d$, let us denote by
$x_{i(.)}$ the ordered statistic computed from the $i^{th}$ coordinate of
$(x_1,\ldots,x_n)$, that is, for all $i=1,\ldots,d$, $$x^i_{i(1)}\leq
x^i_{i(2)}\leq \ldots \leq x^i_{i(n)}.$$ Let $p=\lfloor
n/2\rfloor$ and let 
\begin{eqnarray*}
{\mathcal{B}}_p & = & \left\{\{x_{i(1)};\ldots;x_{i(p)}\} \ ; \ i=1,\ldots,d \
  \mbox{and} \ |\{x_{i(1)};\ldots;x_{i(p)}\}|=p\right\},\\
{\mathcal{B}}_p^c & = & \left\{B\subset  \ \{x_1,\ldots,x_n\}; |B|=p \
  \mbox{and} \ B\notin {\mathcal{B}}_p\right\}.
\end{eqnarray*}
Hence ${\mathcal{B}}_p$ is covered by ${\mathcal{A}}_d$ (by simply taking
$A=\{x^i\leq (x^i_{i(p)}+x^i_{i(p+1)})/2\}$ for each coordinate), and we have that:
$$ |{\mathcal{B}}_p|  \leq  d \ \ \mbox{and} \ \ |{\mathcal{B}}_p^c| \geqslant \binom{n}{p}-d>0.$$
Let $B\in {\mathcal{B}}_p^c$ and $A=\{x^i\leq a\}\in {\mathcal{A}}_d$. If
$|\{x_1,\ldots,x_n\}\cap A|\neq p$, then $\{x_1,\ldots,x_n\}\cap A\neq
B$. Else, since $\{x_1,\ldots,x_n\}\cap A=\{x_j \ ; \ x_j^i\leq a\}$, we have
that $x^i_{i(j)}\leq a$ for all $j=1,\ldots,p$, and $x^i_{i(j)}>a$ for all
$j=p+1,\ldots,n$. So $\{x_1,\ldots,x_n\}\cap A=\{x_{i(1)};\ldots;x_{i(p)}\}$ and
$|\{x_{i(1)};\ldots;x_{i(p)}\}|=p$, leading to $\{x_1,\ldots,x_n\}\cap A\in
{\mathcal{B}}_p$, and then to $\{x_1,\ldots,x_n\}\cap A\neq B$. So, for all $B\in
{\mathcal{B}}_p^c$ and all $A\in {\mathcal{A}}_d$, $\{x_1,\ldots,x_n\}\cap A\neq
B$. \\ So, if
$\displaystyle{\binom{n}{\lfloor n/2\rfloor}> d}$, $(x_1,\ldots,x_n)$ can not be shattered by ${\mathcal{A}}_d$. Thus $$VC({\mathcal{A}}_d)\leq \max\left\{n \ ; \
\binom{n}{\lfloor n/2\rfloor}\leq d\right\}.$$
~
\\

\ni Let $n\geqslant 1$ such that $\displaystyle{\binom{n}{\lfloor n/2\rfloor}\leq
d}$. Let $(x_1,\ldots,x_n)$ be $n$ points of ${\mathbb{R}}^d$ defined as follows: for each coordinate $i=1,\ldots,\binom{n}{\lfloor n/2\rfloor}$, let $\{i_1;\ldots;i_{\lfloor n/2\rfloor}\}$ be the $i^{th}$ subset of $\lfloor n/2\rfloor$ indices in
$\{1;\ldots;n\}$, where the indices are denoted in ascending order, i.e.: $$1\leq i_1<\ldots<i_{\lfloor
  n/2\rfloor}\leq n.$$ Since $\displaystyle{\binom{n}{\lfloor n/2\rfloor}\leq d}$, we obtain $\displaystyle{\binom{n}{\lfloor n/2\rfloor}}$
distinct subsets of indices.\\ Hence we take for each such coordinate $$x_{i_k}^i=k.$$ Then the remaining
values of $(x_1,\ldots,x_n)$ are taken as follows:
\begin{itemize}
\item Since $\displaystyle{\binom{n}{\lfloor n/2\rfloor+1}\leq d}$, for each subset $\{i_1;\ldots;i_{\lfloor n/2\rfloor+1}\}$ of
  $\{1;\ldots;n\}$ with $\lfloor n/2\rfloor+1$ elements, there exists
  $i'\in \{1;\ldots;\binom{n}{\lfloor n/2\rfloor}\}$ such that
  $\{i_1;\ldots;i_{\lfloor n/2\rfloor}\}=\{i'_1;\ldots;i'_{\lfloor
    n/2\rfloor}\}$. Then take $x^{i'}_{i_{\lfloor n/2\rfloor+1}}=\lfloor
  n/2\rfloor+1$. Let us note that, if $n$ is odd, there is a bijection between $i$ and $i'$.
\item Let $\{j_1;\ldots;j_m\}=\{j\notin \{i_1;\ldots;i_{\lfloor n/2\rfloor+1}\}\}$, with $j_1<\ldots<j_m$, and
  let $j_0=i_{\lfloor n/2\rfloor+1}$. Then take $x^{i'}_{j_k}=x^{i'}_{j_{k-1}}+1$.
\end{itemize}
If not filled, the last coordinates are set to be equal to $n$. \\ Hence, we obtain that, for all $j\notin \{i_1;\ldots;i_{\lfloor n/2\rfloor}\}$,
$x_j^i\geqslant \lfloor n/2\rfloor+1$. \\ Then $(x_1,\ldots,x_n)$ is shattered by ${\mathcal{A}}_d$: for $p\in \{0;\ldots;n\}$, let $B=\{x_{i_1};\ldots;x_{i_p}\}\subset  \ \{x_1,\ldots,x_n\}$, with $1\leq i_1<i_2<\ldots <i_p\leq n$ as soon as $p\neq0$.\\ If $p=0$, let $$i_0=\argmin_{1\leq i\leq
  d}\min_j x_j^i,$$ and take $A=\{x^{i_0}\leq \min_j x_j^{i_0}-1\}$. Then
$B=\{x_1,\ldots,x_n\}\cap A=\emptyset$.\\ If $p=n$, let $$i_n=\argmax_{1\leq i\leq
  d}\max_j x_j^i,$$ and take $A=\{x^{i_n}\leq \max_j x_j^{i_n}+1\}$. Then
$B=\{x_1,\ldots,x_n\}\cap A=\{x_1,\ldots,x_n\}$.\\ If $0< p \leq \lfloor
n/2\rfloor$, let $A\in {\mathcal{A}}_d$ be the subset defined by $A=\{x^i\leq
p+1/2\}$, with $i$ the coordinate corresponding to a subset of indices $\{i_1;\dots ; i_{\lfloor n/2\rfloor}\}$ containing $\{i_1; \ldots ; i_p\}$. Then, by definition of $(x^i_1,\ldots,x^i_n)$, 
$B=\{x_1,\ldots,x_n\}\cap A$.\\
If $\lfloor n/2\rfloor+1 \leq p < n$, let $i'$ be the coordinate corresponding to the configuration $\{i_1;\ldots;i_{\lfloor n/2\rfloor+1}\}$ (as defined by $(x_1,\ldots,x_n)$). Let
$A\in {\mathcal{A}}_d$ be the subset defined by $A=\{x^{i'}\leq p+1/2\}$. Then,
by definition of $(x^{i'}_1,\ldots,x^{i'}_n)$, 
$B=\{x_1,\ldots,x_n\}\cap A$.\\ 
Thus $$VC({\mathcal{A}}_d)\geqslant \max\left\{n \ ; \
  \binom{n}{\lfloor n/2\rfloor}\leqslant d\right\}.$$
~
\\

 \ni The bounds are computed thanks to the Stirling's formula: for all $n\geqslant 1$, $$e^{\frac{1}{12n+1}}\sqrt{2\pi n}\left(\frac{n}{e}\right)^n\leq n!\leq e^{\frac{1}{12n}}\sqrt{2\pi n}\left(\frac{n}{e}\right)^n.$$
It follows by a simple computation that, for all $n\geqslant 1$, $\displaystyle{\binom{n}{\lfloor n/2\rfloor}\leqslant \frac{2^{n+1/2}}{\sqrt{\pi}}},$ leading to the lower bound of $VC({\mathcal{A}}_d)$.

\ni Since $VC({\mathcal{A}}_d)\leqslant d$ is increasing with $d$, we have
$VC({\mathcal{A}}_d)\geqslant 3$ for all $d\geqslant 3$. So we will focus only
on integers $3\leqslant n\leqslant d$ to compute the upper bound. We obtain from the Stirling's formula:
\begin{itemize}
\item if $n$ is even
\begin{eqnarray*}
\binom{n}{n/2} & \geqslant & e^{-\frac{9n+1}{6n(12n+1}}\frac{2^{n+1/2}}{\sqrt{\pi n}}\geqslant e^{-\frac{14}{333}}\frac{2^{n+1/2}}{\sqrt{\pi d}},
\end{eqnarray*}
\item if $n$ is odd
\begin{eqnarray*}
 \binom{n}{\lfloor n/2\rfloor} & \geqslant & \left(1-\frac{1}{n^2}\right)^{-\frac{n+1}{2}}\sqrt{\frac{n-1}{n+1}} e^{-\frac{2n}{6(n^2-1)}+\frac{1}{12n+1}}\frac{2^{n+1/2}}{\sqrt{\pi n}}\geqslant  \frac{81}{64\sqrt{2}}e^{-\frac{29}{396}}\frac{2^{n+1/2}}{\sqrt{\pi d}}.
\end{eqnarray*}
\end{itemize}
\ni Thus, it follows that, for all $n\geqslant 3$ such that $\displaystyle{\binom{n}{\lfloor n/2\rfloor} \leq d}$, $\displaystyle{\frac{81}{64\sqrt{2}}e^{-\frac{29}{396}}\frac{2^{n+1/2}}{\sqrt{\pi}}\leq d^{\frac{3}{2}}}$, leading to the upper bound.

%
\end{proof}

\bibliographystyle{apalike}
\bibliography{pruning}

\end{document}